\documentclass[a4paper,12pt]{article}

\usepackage[left=2cm,right=2cm, top=2cm,bottom=2.5cm,bindingoffset=0cm]{geometry}

\usepackage{verbatim}
\usepackage{amsmath}
\usepackage{amsthm}
\usepackage{amssymb}
\usepackage{delarray}
\usepackage{mathrsfs}
\usepackage{graphicx}
\usepackage{cite}

\usepackage{amssymb}

\usepackage{color}

 \DeclareMathOperator{\rank}{rank}  \allowdisplaybreaks

\begin{document}

\begin{center}
{\large\bf STURM--LIOUVILLE-TYPE OPERATORS WITH FROZEN ARGUMENT AND CHEBYSHEV POLYNOMIALS}
\\[0.5cm]

{\large\bf  Tzong-Mo Tsai}
\\[0.1cm]

\footnotesize{General Education Center, Ming Chi University of Technology, New Taipei City, 24301, Taiwan
\\[0.1cm]

\it tsaitm@mail.mcut.edu.tw}
\\[0.2cm]

{\large\bf   Hsiao-Fan Liu}
\\[0.1cm]

\footnotesize{Department of Mathematics, Tamkang University, New Taipei City, 25137, Taiwan
\\[0.1cm]

\it hfliu@mail.tku.edu.tw}
\\[0.2cm]

{\large\bf   Sergey
Buterin}
\\[0.1cm]

\footnotesize{Department of Mathematics, Saratov State University, Astrakhanskaya 83, Saratov 410012, Russia
\\[0.1cm]

\it buterinsa@info.sgu.ru}
\\[0.2cm]

{\large\bf   Lung-Hui Chen}
\\[0.1cm]

\footnotesize{General Education Center, Ming Chi University of Technology, New Taipei City, 24301, Taiwan
\\[0.1cm]

\it mr.lunghuichen@gmail.com}
\\[0.2cm]

{\large\bf  Chung-Tsun Shieh}

\footnotesize{Department of Mathematics, Tamkang University, New Taipei City, 25137, Taiwan
\\[0.1cm]

\it ctshieh@mail.tku.edu.tw}\\[0.2cm]
\end{center}

{\bf Abstract.} The paper deals with Sturm--Liouville-type operators with frozen argument of the form $\ell
y:=-y''(x)+q(x)y(a),$ $y^{(\alpha)}(0)=y^{(\beta)}(1)=0,$ where $\alpha,\beta\in\{0,1\}$ and $a\in[0,1]$ is
an arbitrary fixed rational number. Such nonlocal operators belong to the so-called loaded differential
operators, which often appear in mathematical physics. We focus on the inverse problem of recovering the
potential $q(x)$ from the spectrum of the operator $\ell.$ Our goal is two-fold. Firstly, we establish a deep
connection between the so-called main equation of this inverse problem and Chebyshev polynomials of the first
and the second kinds. This connection gives a new perspective method for solving the inverse problem. In
particular, it allows one to completely describe all non-degenerate and degenerate cases, i.e. when the
solution of the inverse problem is unique or not, respectively. Secondly, we give a complete and convenient
description of iso-spectral potentials in the space of complex-valued integrable functions.


\\

{\large\bf 1. Introduction}
\\

Consider the boundary value problem ${\cal L}:={\cal L}(q(x),a,\alpha,\beta)$ of the form
\begin{equation}\label{1.1}
\ell y:=-y''(x)+q(x)y(a)=\lambda y(x), \quad 0<x<1,
\end{equation}
\begin{equation}\label{1.2}
y^{(\alpha)}(0)=y^{(\beta)}(1)=0,
\end{equation}
where $\lambda$ is the spectral parameter, $q(x)$ is a complex-valued function in $L(0,1),$ to which we refer
as {\it potential}, and $\alpha,\beta\in\{0,1\},$ while $a\in[0,1].$ The operator $\ell$ is called the
Sturm--Liouville-type operator with {\it frozen argument}.

Denote by $\{\lambda_n\}_{n\ge1}$ the spectrum of ${\cal L}$ and consider the following inverse problem.

\medskip
{\bf Inverse Problem 1.} Given $\{\lambda_n\}_{n\ge1},$ $a,$ $\alpha$ and $\beta;$ find $q(x).$

\medskip
Nonlocal operators of the form (\ref{1.1}), (\ref{1.2}) belong to the so-called loaded differential operators
(see, e.g., \cite{Isk, Kral, Nakh76, NakhBor77, Nakh82, Nakh12, Lom14, Lom15, BH21}), which often appear in
mathematical physics.
For example, some models of physical systems with feedback leading to nonlocal differential operators with
frozen argument were described  in \cite{BH21}. The presence of a feedback means that the external affect on
the system depends on its current state. If this state is taken into account only at some fixed physical
point of the system, then mathematically this corresponds to an operator with frozen argument.

Among purely mathematical applications, we illustrate here the so-called {\it method of reduction to loaded
equations} (see, \cite{Nakh82, Nakh12}). For this purpose, let us aim to study the boundary value problem for
the integro-differential equation
\begin{equation}\label{1.2-1}
-y''(x)+\int\limits_0^1H(x,t)y(t)\,\mathrm{d}t=\lambda y(x), \quad 0<x<1,
\end{equation}
subject to boundary conditions (\ref{1.2}) for some $\alpha$ and $\beta.$ The method consists in replacing
equation~(\ref{1.2-1}) with the loaded one
\begin{equation}\label{1.2-2}
-y''(x)+\sum_{\nu=1}^N q_\nu(x)y(a_\nu)=\lambda y(x), \quad 0<x<1,
\end{equation}
possessing frozen arguments $a_1,\ldots,a_N,$ where the sum is an appropriate quadrature formula for
approximating the integral in (\ref{1.2-1}). For example, Simpson's rule (see, e.g., \cite{Atk89})
$$
\int\limits_{x_1}^{x_2} f(x)\,\mathrm{d}x\approx\frac{x_2-x_1}6\Big(f(x_1)+4f\Big(\frac{x_1+x_2}2\Big)
+f(x_2)\Big)
$$
in the case $\alpha=\beta=0$ leads to equation (\ref{1.1}) with $q(x)=2H(x,a)/3$ and $a=1/2.$

Various aspects of Inverse Problem~1 in the case $q(x)\in L_2(0,1)$ were studied in \cite{BBV, BV, BK,
Wang20}. In \cite{BBV, BV, BK}, diverse cases of the triple $(a,\alpha,\beta)$ with rational $a$'s were
considered. In particular, it was established that the inverse problem may be uniquely solvable or not
depending on the parameters $\alpha,\,\beta$ and also on the parity of the integers $k,$ $j$ or $j+k$ taken
from the representation $a=j/k$ under the assumption that $j$ and $k$ are {\it mutually prime}. According to
this, there were highlighted two cases: {\it non-degenerate} and {\it degenerate}~ones, respectively.
Moreover, a complete characterization of the spectrum $\{\lambda_n\}_{n\ge1}$ was given, which includes the
asymptotics for large modulus eigenvalues
along with a special additional condition in the degenerate case.
Specifically, it was established that, in the degenerate case, asymptotically $k$-th part of the spectrum
degenerates, i.e. each $k$-th eigenvalue carries no information on the potential.

For example, in the case when $\alpha=\beta=0$ and $q(x)\in L_2(0,1),$ the spectrum is completely
characterized by the relations
$$
\lambda_n=(\pi n)^2+\varkappa_n, \quad \{\varkappa_n\}\in l_2, \quad \lambda_{kn}=(\pi kn)^2, \quad
n\in{\mathbb N},
$$
i.e. each $k$-th eigenvalue $\lambda_{kn}$ degenerates. In this case, for the unique solvability of Inverse
Problem~1, one should specify the potential on one of the subintervals $((\nu-1)/k,\nu/k),$
$\nu=\overline{1,k}.$ Thus, the smaller part of the spectrum degenerates, the less additional information on
the potential is required, while in the non-degenerate case no extra information is required at all (see
(\ref{deg}) and (\ref{non-deg}) below for a complete description of degenerate and non-degenerate subcases).

This causes instable informativity of the spectrum with respect to~$a,$ which was first revealed in \cite{BV}
(see also \cite{BK}). For example, while a half of the spectrum degenerates as soon as $a=1/2,$ for
$a=a_k:=(k-1)/(2k)$ with even $k$ so does only its $2k$-th part, but $a_k\to1/2$ as $k\to\infty.$

Thus, returning to the method of reduction to loaded equations mentioned above, one can note that this method
is sensitive to choosing the nodes $a_1,\ldots,a_N$ in (\ref{1.2-2}) at least
when approximating the spectrum of the initial integro-differential operator (\ref{1.2-1}) subject to
(\ref{1.2}).

Concerning irrational values of $a,$ it was established in \cite{Wang20} that all they correspond to the
non-degenerate case for all pairs $(\alpha,\beta),$ i.e. the solution of Inverse Problem~1 is always unique
as soon as $a\notin{\mathbb Q}.$ However, the question of the spectrum characterization
still remains open.

In \cite{Nizh-09, AlbHryNizh
} and other works, in connection with the theory of diffusion processes, the case $a=1$ was investigated but
with the special nonlocal boundary conditions
\begin{equation}\label{1.3}
y(0)-\alpha y(1)=y'(1)-\alpha y'(0)+\int\limits_0^1y(t)\overline{q(t)}\,\mathrm{d}t=0, \quad
\alpha\in\{0,1\},
\end{equation}
guarantying the self-adjointness of the corresponding operator generated by (\ref{1.1}) and (\ref{1.3}).
However, such settings never entail the uniqueness of recovering the function $q(x)$ from the spectrum.

In \cite{BH21}, the case of the quasi-periodic boundary conditions of the form
$$
y(0)-\gamma y(1)=y'(0)-\gamma y'(1)=0
$$
for any possible $\gamma\in{\mathbb C}\setminus\{0\}$ was studied, and a complete solution was obtained for
the inverse problem of recovering the potential $q(x)$ from the corresponding spectrum (spectra). Further
aspects of recovering the operator $\ell$ as well as its spectral properties were studied in \cite{XY19-1,
XY19-2, HBY, Hu20}.

In the present paper, we return to Inverse Problem~1. In \cite{BBV, BV, BK}, this problem was reduced to some
linear functional equation with respect to the potential $q(x)$ (see equation~(\ref{dd:1}) in the next
section), which was referred to as {\it  main equation} of the inverse problem. For rational values of~$a,$
the main equation was reduced to linear system (\ref{2.7}) with a special $k\times
k$-matrix~$A_{j,k}^{(\alpha,\beta)},$ whose rank appeared to be ranging between $k-1$ and $k.$ These two
possibilities, in turn, correspond to the degenerate and non-degenerate cases, respectively. In the works
\cite{BBV, BV, BK}, various approaches to studying this matrix and calculating its determinant were
developed.

Here, we establish a deep connection between the matrix $A_{j,k}^{(\alpha,\beta)}$ and Chebyshev polynomials
of the first and the second kinds, which gives another approach for studying Inverse~Problem~1. Using this
new approach, we obtain a complete and convenient description of all iso-spectral complex-valued potentials
in $L(0,1)$ for the degenerate case.

The paper is organized as follows. In the next section, we provide some necessary information on the boundary
value problem ${\cal L}.$ In Section~3, we represent characteristic determinants of the matrices
$A_{1,k}^{(\alpha,\beta)}$ via Chebyshev polynomials. In Section~4, we establish that the matrices
$A_{j,k}^{(\alpha,\beta)}$ can be obtained after substituting $A_{1,k}^{(\alpha,\beta)}$ into appropriate
Chebyshev polynomials, which allows one to study their spectra for $j>0.$ In Section~5, we construct
iso-spectral potentials in the degenerate case. In Section~6, we provide some illustrative examples.
\\

{\large\bf 2. Preliminary information}
\\

Let $C(x,\lambda)$ and $S(x,\lambda)$ be solutions of equation (\ref{1.1}) under the initial conditions
\begin{equation}\label{2.2}
C(a,\lambda)=S'(a,\lambda)=1, \quad S(a,\lambda)=C'(a,\lambda)=0.
\end{equation}
By substitution, it can be easily checked that
\begin{equation}\label{2.1}
C(x,\lambda)=\cos\rho(x-a)+\int\limits_a^x\frac{\sin\rho(x-t)}{\rho}q(t)\,\mathrm{d}t, \quad
S(x,\lambda)=\frac{\sin\rho(x-a)}{\rho}, \quad \rho^2:=\lambda.
\end{equation}
Since these solutions are uniquely determined by conditions~(\ref{2.2}), eigenvalues of the problem~${\cal
L}$ coincide with zeros of the entire function
\begin{equation}\label{2.3}
\Delta_{\alpha,\beta}(\lambda)=
\begin{vmatrix}
C^{(\alpha)}(0,\lambda) & S^{(\alpha)}(0,\lambda) \\
C^{(\beta)}(1,\lambda) & S^{(\beta)}(1,\lambda)
\end{vmatrix},
\end{equation}
which is called {\it characteristics function} of ${\cal L}.$
%

Without loss of generality, we always assume that $0\le a\le1/2$ since the spectrum of the problem ${\cal
L}(q(x),a,\alpha,\beta),$ obviously, coincides with the one of ${\cal L}(q(1-x),1-a,\beta,\alpha).$

Substituting (\ref{2.1}) into (\ref{2.3}) one can obtain the following representations (see, e.g.,
\cite{BK}):
\begin{equation} \label{2.5}
\Delta_{\alpha,\alpha}(\lambda)=\rho^{2\alpha}\Big(\frac{\sin\rho}{\rho}+
\int\limits_0^1W_{\alpha,\alpha}(x)\frac{\cos\rho x}{\rho^2}\,\mathrm{d}x\Big),\;\;
\Delta_{\alpha,\beta}(\lambda)=(-1)^\alpha\cos\rho+\int\limits_0^1W_{\alpha,\beta}(x)\frac{\sin\rho
x}{\rho}\,\mathrm{d}x
\end{equation}
for $\alpha\ne\beta,$ where the functions $W_{\alpha,\beta}(x)$ have the form
\begin{equation}\label{dd:1}
\displaystyle W_{\alpha,\beta}(x)=\frac{(-1)^{\alpha\beta}}{2} \left\{
\begin{array}{cl}
\displaystyle  q(1-a+x)+dq(1-a-x), & x\in(0,a),\\[3mm]
\displaystyle cq(1+a-x)+dq(1-a-x), & x\in(a,1-a),\\[3mm]
\displaystyle c\Big(q(1+a-x)+q(x-1+a)\Big), & x\in(1-a,1),
\end{array}\right.
\end{equation}
while the numbers $c$ and $d$ are determined by the formulae
\begin{equation}\label{2.8.1}
c = (-1)^{\beta + 1}, \quad d=(-1)^{\alpha+\beta}.
\end{equation}
Assuming the function $W_{\alpha,\beta}(x)$ to be known, one can consider (\ref{dd:1}) as a linear functional
equation with respect to $q(x).$ Since each function $\Delta_{\alpha,\beta}(\lambda)$ is uniquely determined
by its zeros:
\begin{equation}\label{EQ}
\Delta_{\alpha,\beta}(\lambda) =(-1)^\alpha(\lambda_1-\lambda)^{\alpha\beta}
\prod\limits_{n=1+\alpha\beta}^\infty\frac{\lambda_n-\lambda}{\Big(n-\frac{\alpha+\beta}2\Big)^2\pi^2},
\end{equation}
Inverse Problem~1 is equivalent to this functional equation (\ref{dd:1}), which is called {\it main equation}
of the inverse problem.

If $a$ is rational, i.e. there exist {\it mutually prime} integers $j$ and $k$ such that $a=j/k,$ then the
main equation can be represented in the following way:
\begin{equation}\label{2.7}
W_{\alpha,\beta}(x)=\frac{(-1)^{\alpha\beta}}{2} Q^{-1}A_{j,k}^{(\alpha,\beta)} Rq(x), \quad 0<x<1,
\end{equation}
where $A_{j,k}^{(\alpha,\beta)}$ is a square matrix of order $k,$ while $Q$ and $R$ are bijective operators
mapping $L(0,1)$ onto $(L(0,b))^k,$ $b:=1/k,$ and acting by the formulae
\begin{equation} \label{2.4.1}
Qf:=(f,Q_2f,\ldots,Q_kf)^T, \quad Rf:=(R_1f,R_2f,\ldots,R_kf)^T.
\end{equation}
Here, $T$ is the transposition sign, while $Q_\nu$ and $R_\nu$ are shift and involution operators mapping
$L(0,1)$ onto $L(0,b),$ which are determined by the formulae
\begin{equation}\label{2.3.1}
Q_\nu f(x)=\left\{
\begin{array}{ll}
\displaystyle\!\! f((\nu-1)b+x) \;\; \text{for odd} \; \nu,\\[3mm]
\displaystyle\!\! f(\nu b-x) \;\; \text{for even} \; \nu,
\end{array}\right.
R_\nu f(x)=\left\{
\begin{array}{l}
\displaystyle\!\! f((k-\nu)b+x) \;\; \text{for even} \; j+\nu,\\[3mm]
\displaystyle\!\! f((k-\nu+1)b-x) \;\; \text{for odd} \; j+\nu,
\end{array}\right.
\!\!
\end{equation}
where $x\in(0,b)$ and $\nu=\overline{1,k}.$

The matrix $A_{j,k}^{(\alpha,\beta)}=(a_{m,n})_{m,n=\overline{1,k}},$ in turn, is constructed in the
following way. For $k=1,$ it consists of a single element $a_{1,1}=2(-1)^{\beta+1}\alpha$ as soon as $j=0,$
while, for $k\ge2j\ge2,$ its elements are determined by the formulae
\begin{equation}\label{2.8}
\left.\begin{array}{rll}
(i) & a_{m,j-m+1}= 1, & m=\overline{1,j},\\[3mm]
(ii) & a_{m,m+j}= d, & m = \overline{1,k-j}, \\[3mm]
(iii) & a_{m,m-j}= c, & m=\overline{j +1,k},\\[3mm]
(iv) & a_{m,2k-m-j+1}= c, & m=\overline{k-j + 1,k},\\[3mm]
(v) & a_{m,n}=0 & \text{for the remaining pairs}\;\; (m,n).
\end{array}\right\}
\end{equation}
The items $(i)$--$(iv)$ in \eqref{2.8} correspond to subdiagonals, consisting of equal elements: $1,$ $c$
and~$d.$ For example, the the matrix $A_{3,7}^{(\alpha,\beta)}$ has the form
\begin{equation*}
\begin{pmatrix}
\cdot & \cdot & 1 & d & \cdot & \cdot & \cdot \\
\cdot & 1 & \cdot & \cdot & d & \cdot & \cdot \\
1 & \cdot & \cdot & \cdot & \cdot & d & \cdot \\
c & \cdot & \cdot & \cdot & \cdot & \cdot & d \\
\cdot & c & \cdot & \cdot & \cdot & \cdot & c \\
\cdot & \cdot &c & \cdot & \cdot & c & \cdot \\
\cdot & \cdot & \cdot & c & c & \cdot & \cdot
\end{pmatrix},
\end{equation*}
where dots indicate zero elements.

In papers \cite{BBV, BV, BK}, various approaches were used for calculating the determinant and the rank of
$A_{j,k}^{(\alpha,\beta)},$ depending on generality of the situation. In particular, in \cite{BV}, a
reduction-type algorithm was suggested for the case $\alpha=\beta=0.$ This algorithm appeared to be
equivalent to the Euclidean algorithm for finding the greatest common devisor of the numbers $j$ and $k.$
Since $j$ and $k$ are mutually prime, the algorithm gave consecutive relations leading to the result:
$$
\rank A_{j,k}^{(0,0)}=\rank A_{j,k-j}^{(0,0)}+j= \rank A_{k-2j,k-j}^{(0,0)}+j=\ldots=\rank
A_{0,1}^{(0,0)}+k-1 =k-1.
$$
Later, in \cite{BK}, it was established that the rank of $A_{j,k}^{(\alpha,\beta)}$ cannot be less than
$k-1.$ For this purpose, a combinatorial approach for calculating the determinant of
$A_{j,k}^{(\alpha,\beta)}$ was suggested, which was based on studying properties of an undirected graph $G$
corresponding to a special traversal of nonzero elements of $A_{j,k}^{(\alpha,\beta)}.$ It was established
that $G$ was a bipartite Eulerian cycle, which has led to representing $\det A_{j,k}^{(\alpha,\beta)}$ as a
sum of precisely two products consisting of nonzero elements of $A_{j,k}^{(\alpha,\beta)}.$ This gave a
complete classification of degenerate and non-degenerate cases corresponding to non-unique and unique
solvability of Inverse Problem~1, respectively. Specifically, the degenerate case occurs when one of the
following groups of conditions is fulfilled:
\begin{equation}\label{deg}
\left.\begin{array}{rl}
\text{(I)}& \alpha=\beta=0;\\[3mm]
\text{(II)}& \alpha=0, \ \beta=1 \text{ and } j \text{ is even;}\\[3mm]
\text{(III)}& \alpha=1,\ \beta=0 \text{ and } k+j \text{ is even;}\\[3mm]
\text{(IV)}& \alpha=\beta=1 \text{ and } k \text{ is even};
\end{array}\right\}
\end{equation}
while the non-degenerate case includes the remaining groups of conditions:
\begin{equation}\label{non-deg}
\left.\begin{array}{rl}
\text{(V)}& \alpha=0, \;\beta=1 \text{ and } j \text{ is odd;}\\[3mm]
\text{(VI)}& \alpha=1,\ \beta=0 \text{ and } k+j \text{ is odd;}\\[3mm]
\text{(VII)}& \alpha=\beta=1 \text{ and } k \text{ is odd.}
\end{array}\right\}
\end{equation}
This classification remains valid also for $a>1/2,$ i.e. it holds for all relevant $j\in\{0,\ldots,k\}.$

In the subsequent sections, we establish a deep connection between the matrix $A_{j,k}^{(\alpha,\beta)}$ and
Chebyshev polynomials of the first and the second kinds. This connection gives, in particular, another
approach for studying the main equation~(\ref{2.7}). Using this approach one can easily give a complete
description of iso-spectral potentials in the degenerate case (see Section~5).
\\

{\large\bf 3. Chebyshev polynomials and the case $j=1$}
\\

First, we give some necessary information about Chebyshev polynomials $T_{n}(z)$ and $U_{n}(z)$ of the first
and the second kinds, respectively, which can be defined by the formulae
\begin{equation}\label{3.10,1}
T_{n}(\cos\theta)=\cos n\theta, \quad U_{n}(\cos\theta)=\frac{\sin(n+1)\theta}{\sin\theta}, \quad
n+1\in{\mathbb N}.
\end{equation}
Alternatively, one can use the following recurrent relation
\begin{equation}\label{3.8,1}
Y_{n+1}(z)=2z Y_{n}(z)- Y_{n-1}(z).
\end{equation}
Then the polynomials of the first kind $T_n(z)=Y_n(z)$ are determined by the initial conditions
\begin{equation}\label{3.8,2}
T_{0}(z)=1, \quad T_{1}(z)=z,
\end{equation}
while the initial conditions
\begin{equation}\label{3.9,1}
U_{0}(z)=1, \quad U_{1}(z)=2z
\end{equation}
determine the polynomials of the second kind  $Y_n(z)=U_n(z).$ For more details, see, e.g., \cite{Riv90}.

It is well known and also can be easily seen that Chebyshev polynomials may possess only simple zeros, and
they are always odd or even functions in accordance with the parity of $n.$ In particular, we have
$T_{n}(0)=U_{n}(0)=0$ as soon as $n$ is odd, and $T_{n}(0)U_{n}(0)\ne0$ for even $n.$

Let us proceed with studying the matrix $A_{j,k}^{(\alpha,\beta)}$ for $k\ge2.$ In this section, we focus on
the case~$j=1.$ Consider the characteristic polynomial
\begin{equation}\label{2.4}
p_k(z):=\det(zI-A_{1,k}^{(\alpha,\beta)}),
\end{equation}
where $I$ is the unit matrix. The following lemma holds.

\medskip
{\bf Lemma 1. }{\it The characteristic polynomial of the matrix $A_{1,k}^{(\alpha,\beta)}$ has the form
\begin{equation}\label{2.3.2}
p_k(z)=(z-c)q_{k-1}(z)-cdq_{k-2}(z),
\end{equation}
where $c$ and $d$ are determined by (\ref{2.8.1}), while the polynomials $q_n(z)$ can be found from the
recurrent relations
\begin{equation}\label{2.3.3}
q_0(z)=1, \quad q_1(z)=z-1, \quad q_{n+1}(z)=z q_n(z)- cdq_{n-1}(z), \quad n=\overline{1,k-2}.
\end{equation}
}

{\it Proof.} First, we note that $q_\nu(z)$ is the characteristic polynomial of the three-diagonal matrix
$B_\nu$ that is obtained from $A_{1,\nu+1}^{(\alpha,\beta)}$ by removing the last column along with the last
row, i.e.
\begin{equation}\label{2.3.4}
q_\nu(z)=\det(zI-B_\nu)= \left|\begin{array}{ccccc}
z-1& -d &        &        &    \\[1mm]
-c & z  &   -d   &        &    \\[1mm]
   & -c & \ddots & \ddots &    \\[1mm]
   &    & \ddots &   z    & -d \\[1mm]
   &    &        &  -c    & z  \\[1mm]
\end{array}\right|,
\end{equation}
where each of both subdiagonals consists of equal elements, while all elements of the main diagonal starting
from the second position are equal too. Indeed, for $\nu=1$ formula (\ref{2.3.4}) is obvious. Further, let it
hold for any $\nu\le n.$ Then expanding the determinant in (\ref{2.3.4}) for $\nu=n+1$ with respect to the
elements of the last row we obtain the last equality in (\ref{2.3.3}).

Finally, expanding the determinant
$$
\det(zI-A_{1,k}^{(\alpha,\beta)})= \left|\begin{array}{ccccc}
z-1& -d &        &        &    \\[1mm]
-c & z  &   -d   &        &    \\[1mm]
   & -c & \ddots & \ddots &    \\[1mm]
   &    & \ddots &   z    & -d \\[1mm]
   &    &        &  -c    & z-c  \\[1mm]
\end{array}\right|
$$
with respect to the last row, we obtain representation (\ref{2.3.2}). $\hfill\Box$

\medskip
The following corollary gives the classifications (\ref{deg}) and (\ref{non-deg}) for $j=1.$

\medskip
{\bf Corollary 1. }{\it The determinant of the matrix $A_{1,k}^{(\alpha,\beta)}$ can be calculated by the
formula}
\begin{equation}\label{3.5}
\det A_{1,k}^{(\alpha,\beta)} = \left\{
\begin{array}{l}
\displaystyle (-cd)^{(k-1)/2}(1+c) \quad \text{if}\;\;k\;\;\text{is odd},\\[1mm]
\displaystyle c(-cd)^{k/2-1}(1-d) \quad \text{if}\;\;k\;\;\text{is even}.
\end{array}\right.
\end{equation}

{\it Proof.} According to (\ref{2.4}) and (\ref{2.3.2}), we have $\det A_{1,k}^{(\alpha,\beta)}=(-1)^{k+1}c
(q_{k-1}(0)+dq_{k-2}(0)).$ The first two formulae in (\ref{2.3.3}) give $q_0(0)=1$ and $q_1(0)=-1.$ Assume
that
\begin{equation}\label{3.6}
q_{2\nu}(0)=(-cd)^\nu, \quad q_{2\nu+1}(0)=-(-cd)^\nu, \quad 0\le\nu\le l,
\end{equation}
for some $l\in{\mathbb N}.$ Then the last relation in (\ref{2.3.3}) implies
$$
q_{2(l+1)}(0)=-cdq_{2l}(0)=(-cd)^{l+1}, \quad q_{2(l+1)+1}(0)=-cdq_{2l+1}(0)=-(-cd)^{l+1}.
$$
Hence, (\ref{3.6}) holds for all $\nu\ge0.$ Substituting (\ref{3.6}) into the first formula of this proof, we
arrive~at
$$
\det A_{1,k}^{(\alpha,\beta)} = \left\{
\begin{array}{l}
\displaystyle c(q_{2\nu}(0)+dq_{2\nu-1}(0))=(-cd)^\nu(1+c) \quad \text{for}\;\;k=2\nu+1,\\[2mm]
\displaystyle -c(q_{2\nu+1}(0)+dq_{2\nu}(0))=c(-cd)^\nu(1-d) \quad \text{for}\;\;k=2\nu+2,
\end{array}\right.
$$
which finalizes the proof. $\hfill\Box$

\medskip
The main result of this section is contained in the following theorem.

\medskip
{\bf Theorem 1. }{\it The following representations hold:}
\begin{equation}\label{3.6.1}
\det(zI-A_{1,k}^{(0,0)})=i^{k-1}zU_{k-1}\Big(\frac{z}{2i}\Big),
\end{equation}

\begin{equation}\label{3.7}
\det(zI-A_{1,k}^{(0,1)})=2i^kT_k\Big(\frac{z}{2i}\Big)-2i^{k-1}U_{k-1}\Big(\frac{z}{2i}\Big),
\end{equation}

\begin{equation}\label{3.8}
\det(zI-A_{1,k}^{(1,0)})=2T_k\Big(\frac{z}2\Big),
\end{equation}

\begin{equation}\label{3.9}
\det(zI-A_{1,k}^{(1,1)})=(z-2)U_{k-1}\Big(\frac{z}2\Big).
\end{equation}

{\it Proof.} First, we note that, for $Y_n(z)=T_n(z)$ and $Y_n(z)=U_n(z),$ relations
(\ref{3.8,1})--(\ref{3.9,1}) give
\begin{equation}\label{3.14}
2zU_n(z)=\left\{\begin{array}{r} 2z, \quad n=0,\\[2mm]
4z^2,\quad n=1,
\end{array}\right. \quad
2T_{n+1}(z)+2iU_n(z)=\left\{\begin{array}{r} 2z+2i, \quad n=0,\\[2mm]
4z^2+4iz-2,\quad n=1,
\end{array}\right.
\end{equation}
\begin{equation}\label{3.15}
2T_{n+1}(z)=\left\{\begin{array}{r} 2z, \quad n=0,\\[2mm]
4z^2-2,\quad n=1,
\end{array}\right. \quad
(2z-2)U_n(z)=\left\{\begin{array}{r} 2z-2, \quad n=0,\\[2mm]
4z^2-4z,\quad n=1.
\end{array}\right.
\end{equation}

Let $\alpha=0.$ Then formulae (\ref{2.8.1}) give $cd=-1.$ Then (\ref{2.4}) and (\ref{2.3.2}) imply
\begin{equation}\label{3.10}
\det(zI-A_{1,k}^{(0,\beta)})=(z-c)q_{k-1}(z)+q_{k-2}(z).
\end{equation}
Put $Y_n(z):=i^{-n}q_n(2iz),$ $n=\overline{0,k-1}.$ Using (\ref{2.3.3}), one can easily check that the
polynomials $Y_n(z)$ satisfy the recurrent relations~(\ref{3.8,1}). Substituting $q_n(z)=i^nY_n(z/(2i))$ into
(\ref{3.10}), we get
\begin{equation}\label{3.12}
\det(zI-A_{1,k}^{(0,\beta)})=i^k\Big(\Big(2\frac{z}{2i}+ ic\Big)Y_{k-1}\Big(\frac{z}{2i}\Big)
-Y_{k-2}\Big(\frac{z}{2i}\Big)\Big).
\end{equation}
Using (\ref{2.3.3}), we calculate: $Y_{-1}(z)=-i,$ $Y_0(z)=1$ and $Y_1(z)=2z+i.$ Hence, we obtain
$$
(2z+ic)Y_n(z)-Y_{n-1}(z)=\left\{\begin{array}{r} \left.\begin{array}{r} 2z, \quad n=0,\\[2mm]
4z^2,\quad n=1,
\end{array}\right\} \quad \beta=0,\\[5mm]
\left.\begin{array}{r} 2z+2i, \quad n=0,\\[2mm]
4z^2+4iz-2,\quad n=1,
\end{array}\right\} \quad \beta=1.
\end{array}\right.
$$
Comparing this with (\ref{3.14}), we get
$$
(2z+ic)Y_n(z)-Y_{n-1}(z)=\left\{\begin{array}{r} 2zU_n(z), \quad \beta=0,\\[5mm]
2T_{n+1}(z)+2iU_n(z), \quad \beta=1,
\end{array}\right. \quad n=\overline{0,k-1},
$$
which along with (\ref{3.12}) gives (\ref{3.6.1}) and (\ref{3.7}).

Further, let $\alpha=1.$  Then $cd=1,$ and formulae (\ref{2.4}) and (\ref{2.3.2}) imply
\begin{equation}\label{2.10.1}
\det(zI-A_{1,k}^{(1,\beta)})=(z-c)q_{k-1}(z)-q_{k-2}(z).
\end{equation}
Put $Y_n(z):=q_n(2z),$ $n=\overline{0,k-1}.$ By virtue of (\ref{2.3.3}), these polynomials $Y_n(z)$ satisfy
the recurrent relations~(\ref{3.8,1}). Substituting $q_n(z)=Y_n(z/2)$ into (\ref{2.10.1}), we get
\begin{equation}\label{2.12.1}
\det(zI-A_{1,k}^{(1,\beta)})=\Big(2\frac{z}2 -c\Big)Y_{k-1}\Big(\frac{z}2\Big) -Y_{k-2}\Big(\frac{z}2\Big).
\end{equation}
By virtue of (\ref{2.3.3}), we have
$$
(2z-c)Y_n(z)-Y_{n-1}(z)=\left\{\begin{array}{r} \left.\begin{array}{r} 2z, \quad n=0,\\[2mm]
4z^2-2,\quad n=1,
\end{array}\right\} \quad \beta=0,\\[5mm]
\left.\begin{array}{r} 2z-2, \quad n=0,\\[2mm]
4z^2-4z,\quad n=1,
\end{array}\right\} \quad \beta=1.
\end{array}\right.
$$
Comparing this with (\ref{3.15}) and using (\ref{2.12.1}), we arrive at (\ref{3.8}) and (\ref{3.9}).
$\hfill\Box$

\medskip
{\bf Corollary 2. }{\it Denote by $\sigma(A)$ the spectrum of the matrix $A.$ Then
\begin{equation}\label{3.16}
\sigma(A_{1,k}^{(0,0)})=\{0\}\cup\Big\{2i\cos\frac{\nu\pi}{k}\Big\}_{\nu=\overline{1,k-1}},
\end{equation}
\begin{equation}\label{3.17}
0\notin \sigma(A_{1,k}^{(0,1)}),
\end{equation}
\begin{equation}\label{3.18}
\sigma(A_{1,k}^{(1,0)})=\Big\{2\cos\frac{(2\nu+1)\pi}{2k}\Big\}_{\nu=\overline{0,k-1}},
\end{equation}
\begin{equation}\label{3.19}
\sigma(A_{1,k}^{(1,1)})=\Big\{2\cos\frac{\nu\pi}{k}\Big\}_{\nu=\overline{0,k-1}}.
\end{equation}
}

{\it Proof.} It is well known and also can be obtained as a simple corollary from (\ref{3.10,1}) that the
sets of zeros of the polynomials $T_n(z)$ and $U_n(z)$ have the forms
\begin{equation}\label{3.19-1}
{\cal T}_n:=
\Big\{\cos\frac{(2\nu+1)\pi}{2n} \Big\}_{\nu=\overline{0,n-1}}, \quad {\cal
U}_n:=
\Big\{\cos\frac{\nu\pi}{n+1} \Big\}_{\nu=\overline{1,n}},
\end{equation}
respectively. Thus, (\ref{3.16}), (\ref{3.18}) and (\ref{3.19}) follow directly from (\ref{3.6.1}),
(\ref{3.8}) and (\ref{3.9}). Concerning~(\ref{3.17}), it is sufficient to recall that, in (\ref{3.7}),
$T_k(0)U_{k-1}(0)=0,$ while $T_k(0)\ne U_{k-1}(0).$ $\hfill\Box$
\\

{\large\bf 4. The case $j>1$}
\\

In this section, we establish connections between the matrices $A_{j,k}^{(\alpha,\beta)}$ and
$A_{1,k}^{(\alpha,\beta)},$ which allow one to reduce studying the case $j>1$ to the case $j=1.$ Namely, the
following theorem holds.

\medskip
{\bf Theorem 2. }{\it For $\beta=0,1$ and $j=\overline{1,n_k},$ where $n_k=[k/2],$ the following relations
hold:
\begin{equation}\label{4.1}
A_{j,k}^{(0,\beta)}=U_{j-1}\Big(-\frac{c}2A_{1,k}^{(1,1-\beta)}\Big)A_{1,k}^{(0,\beta)},
\end{equation}
\begin{equation}\label{4.2}
A_{j,k}^{(1,\beta)}=2c T_j\Big(\frac{c}2A_{1,k}^{(1,\beta)}\Big).
\end{equation}
Here, $[\,x\,]$ denotes the integer part of $x$ and, as before, $c=(-1)^{1+\beta}.$ }

\medskip
{\it Proof.} For $j=1,$ the assertion is obvious. According to the formulae $U_1(z)=2z$ and $T_2(z)=2z^2-1,$
each of relations (\ref{4.1}) and (\ref{4.2}) for $j=2$ is equivalent to the common relation
\begin{equation}\label{4.3}
A_{2,k}^{(\alpha,\beta)}=dA_{1,k}^{(1,\gamma)}A_{1,k}^{(\alpha,\beta)} -2\alpha cI, \quad
\gamma:=\left\{\begin{array}{r}1-\beta, \;\alpha=0,\\[1mm]
\beta,\;\alpha=1,
\end{array}\right. \quad
\alpha, \beta=0,1.
\end{equation}
Consider a column vector $X=(x_1,\ldots,x_k)^T$ and denote $[X]_n:=x_n$ for $n=\overline{1,k}.$ Then, by
virtue of (\ref{2.8}), we have
\begin{equation}\label{4.6}
[A_{j,k}^{(\alpha,\beta)}X]_m=\left\{\begin{array}{cl}
x_{j-m+1} +dx_{j+m}, & m=\overline{1,j},\\[2mm]
cx_{m-j} +dx_{j+m}, & m=\overline{j+1,k-j},\\[2mm]
c(x_{m-j} +x_{2k-m-j+1}), & m=\overline{k-j+1,k}.
\end{array}\right.
\end{equation}
In particular, this gives the formulae
\begin{equation}\label{4.6.1}
[A_{1,k}^{(\alpha,\beta)}X]_m=\left\{\begin{array}{cl}
x_1 +dx_2, & m=1,\\[2mm]
cx_{m-1} +dx_{m+1}, & m=\overline{2,k-1},\\[2mm]
c(x_{k-1} +x_k), & m=k,
\end{array}\right.
\end{equation}
\begin{equation}\label{4.9}
[A_{1,k}^{(1,\gamma)}X]_m=\left\{\begin{array}{cl}
x_1 +dx_2, & m=1,\\[2mm]
d(x_{m-1} +x_{m+1}), & m=\overline{2,k-1},\\[2mm]
d(x_{k-1} +x_k), & m=k,
\end{array}\right.
\end{equation}
since $\gamma=\alpha\beta+(1-\alpha)(1-\beta)=2\alpha\beta-\alpha-\beta+1$ and, hence,
$(-1)^{1+\gamma}=(-1)^{\alpha+\beta}=d.$ Substituting $A_{1,k}^{(\alpha,\beta)}X$ given by (\ref{4.6.1}) into
(\ref{4.9}) instead of $X,$ we get the relation
$$
[A_{1,k}^{(1,\gamma)}A_{1,k}^{(\alpha,\beta)}X]_m=\left\{\begin{array}{ll}
(1+cd)x_m +dx_{3-m}+x_{2+m}, & m=1,2,\\[2mm]
(1+cd)x_m +cdx_{m-2} +x_{m+2}, & m=\overline{3,k-2},\\[2mm]
(1+cd)x_m +cd(x_{m-2} +x_{2k-m-1}), & m=k-1,k.
\end{array}\right.
$$
Comparing this with (\ref{4.6}) for $j=2$ and taking into account that $c+d=2\alpha c,$ we arrive at
(\ref{4.3}).

Assume now that (\ref{4.1}) and (\ref{4.2}) are valid when $j=\overline{1,\nu}$ for some $\nu\in\{2,\ldots,
n_k-1\}.$ Then, according to (\ref{3.8,1}), relation (\ref{4.1}) for $j=\nu+1$ is equivalent to
$$
A_{\nu+1,k}^{(0,\beta)}=-cA_{1,k}^{(1,1-\beta)}U_{\nu-1}\Big(-\frac{c}2A_{1,k}^{(1,1-\beta)}\Big)
A_{1,k}^{(0,\beta)} -U_{\nu-2}\Big(-\frac{c}2A_{1,k}^{(1,1-\beta)}\Big) A_{1,k}^{(0,\beta)}
\qquad\qquad\qquad
$$
\begin{equation}\label{4.10}
\qquad\qquad\qquad\qquad\qquad\qquad\qquad\qquad\qquad\qquad\qquad\qquad\quad\;\,
=-cA_{1,k}^{(1,1-\beta)}A_{\nu,k}^{(0,\beta)} -A_{\nu-1,k}^{(0,\beta)},
\end{equation}
while (\ref{4.2}) for $j=\nu+1$ takes the form
\begin{equation}\label{4.11}
A_{\nu+1,k}^{(1,\beta)}=2A_{1,k}^{(1,\beta)}T_\nu\Big(\frac{c}2A_{1,k}^{(1,\beta)}\Big) -2c
T_{\nu-1}\Big(\frac{c}2A_{1,k}^{(1,\beta)}\Big) =cA_{1,k}^{(1,\beta)}A_{\nu,k}^{(1,\beta)}
-A_{\nu-1,k}^{(1,\beta)}.
\end{equation}
Using the relation $(-1)^{\alpha+1} c=d$ along with the definition of $\gamma$ in~(\ref{4.3}), one can
rewrite~(\ref{4.10}) and~(\ref{4.11}) in the following common form:
\begin{equation}\label{4.12}
A_{\nu+1,k}^{(\alpha,\beta)}=dA_{1,k}^{(1,\gamma)}A_{\nu,k}^{(\alpha,\beta)} -A_{\nu-1,k}^{(\alpha,\beta)},
\quad \alpha,\beta=0,1.
\end{equation}
Thus, it remains to prove relation (\ref{4.12}).  Using (\ref{4.6}), we calculate
$$
[(A_{\nu+1,k}^{(\alpha,\beta)}+A_{\nu-1,k}^{(\alpha,\beta)})X]_m
\qquad\qquad\qquad\qquad\qquad\qquad\qquad\qquad\qquad\qquad\qquad\qquad\qquad\qquad
$$
\begin{equation}\label{4.13}
\qquad\qquad\quad\;\,=\left\{\begin{array}{cl}
x_{\nu-m+2} +dx_{\nu+1+m}+x_{\nu-m} +dx_{\nu-1+m}, & m=\overline{1,\nu-1},\\[2mm]
x_{\nu-m+2} +dx_{\nu+1+m}+cx_{m-\nu+1} +dx_{\nu-1+m}, & m=\nu,\nu+1,\\[2mm]
cx_{m-\nu-1} +dx_{\nu+1+m}+cx_{m-\nu+1} +dx_{\nu-1+m}, & m=\overline{\nu+2,k-\nu-1},\\[2mm]
c(x_{m-\nu-1} +x_{2k-m-\nu})+cx_{m-\nu+1} +dx_{\nu-1+m}, & m=k-\nu,k-\nu+1,\\[2mm]
c(x_{m-\nu-1} +x_{2k-m-\nu}+x_{m-\nu+1} +x_{2k-m-\nu+2}), & m=\overline{k-\nu+2,k}.
\end{array}\right.
\end{equation}
Further, substituting $A_{\nu,k}^{(\alpha,\beta)}X$ given by (\ref{4.6}) into (\ref{4.9}) instead of $X,$ we
get the relation
$$
[A_{1,k}^{(1,\gamma)}A_{\nu,k}^{(\alpha,\beta)}X]_m=\left\{\begin{array}{cl}
d(x_{\nu-m+2} +dx_{\nu-1+m}+x_{\nu-m} +dx_{\nu+1+m}), & m=\overline{1,\nu-1},\\[2mm]
d(x_{\nu-m+2} +dx_{\nu-1+m}+cx_{m-\nu+1} +dx_{\nu+1+m}), & m=\nu,\nu+1,\\[2mm]
d(cx_{m-\nu-1} +dx_{\nu-1+m}+cx_{m-\nu+1} +dx_{\nu+1+m}), & m=\overline{\nu+2,k-\nu-1},\\[2mm]
d(cx_{m-\nu-1} +dx_{\nu-1+m}+c(x_{m-\nu+1} +x_{2k-m-\nu})), & m=k-\nu,k-\nu+1,\\[2mm]
dc(x_{m-\nu-1} +x_{2k-m-\nu+2}+x_{m-\nu+1} +x_{2k-m-\nu}), & m=\overline{k-\nu+2,k}.
\end{array}\right.
$$
Comparing this with (\ref{4.13}), we arrive at (\ref{4.12}). $\hfill\Box$

\medskip
{\bf Corollary 3. }{\it For $1\le j\le[k/2],$ the matrix $A_{j,k}^{(\alpha,\beta)}$ is degenerate, i.e. $\det
A_{j,k}^{(\alpha,\beta)}=0,$ if and only if one of the four conditions in (\ref{deg}) is fulfilled.

Equivalently, $\det A_{j,k}^{(\alpha,\beta)}\ne0$ if and only if one of the three conditions in
(\ref{non-deg}) is fulfilled.}

\medskip
{\it Proof.} Consider $\alpha=\beta=0$ first. Then (\ref{3.16}) and (\ref{4.1})
imply $\det A_{j,k}^{(0,0)}=0$ for any possible~$j$ and~$k.$ Thus, the assertion of the corollary is proven
for condition (I) in (\ref{deg}).

The rest part of the proof is based on the following well-known fact, which is valid both for Hermitian and
non-Hermitian square matrices $A,$ being a particular case of the corresponding abstract assertion (see,
e.g., Theorem 3.3 on p.~16 in \cite{GGK90} or Theorem 10.28 on p.~263 in \cite{Rud91}).

\medskip
{\bf Proposition 1. }{\it Let $P(z)$ be an algebraic polynomial and $A$ be a square matrix. Then
$$
\sigma(P(A))=\{P(z)\}_{z\in\sigma(A)}.
$$

Moreover, if $X$ is an eigenvector corresponding to an eigenvalue $z_0$ of the matrix $A,$ then~$X$ is an
eigenvector related to the eigenvalue $P(z_0)$ of $P(A).$


}

\medskip
Let us return to the proof of Corollary~3. For $\alpha=0$ and $\beta=1,$ according to (\ref{2.8.1}),
(\ref{3.17}) and (\ref{4.1}),
we have $\det A_{j,k}^{(0,1)}=0$ if and only if $0\in\sigma(U_{j-1}((-1/2)A_{1,k}^{(1,0)})).$ By virtue of
Proposition~1, this inclusion is equivalent to the relation ${\cal
U}_{j-1}\cap\sigma((-1/2)A_{1,k}^{(1,0)})\ne\emptyset,$ where, as in the proof of Corollary~2, we use the
designation ${\cal U}_n=\{z:U_n(z)=0\}$ similarly to ${\cal T}_n=\{z:T_n(z)=0\}.$ Thus, according to
(\ref{3.18}) and (\ref{3.19-1}), the latter intersection is not empty if and only if
$$
\cos\frac{\nu\pi}j+\cos\frac{(2l+1)\pi}{2k}=0
$$
for a certain choice of $\nu\in\{1,\ldots,j-1\}$ and $l\in\{0,\ldots,k-1\}.$ The latter, in turn, is
equivalent to the relation
\begin{equation}\label{4.14}
\frac\nu{j}+(-1)^s\frac{2l+1}{2k}=1+2m
\end{equation}
for some integers $s$ and $m.$ Obviously, (\ref{4.14}) implies the evenness of $j.$ Conversely, let $j$ be
even. Then $k$ is odd, and (\ref{4.14}) holds for $s=m=0$ as soon as $\nu=j/2$ and $l=(k-1)/2.$

For $\alpha=1$ and $\beta=0,$ relations (\ref{2.8.1}) and (\ref{4.2}) imply that $\det A_{j,k}^{(1,0)}=0$ is
equivalent to ${\cal T}_j\cap\sigma((-1/2)A_{1,k}^{(1,0)})\ne\emptyset.$ By virtue of (\ref{3.18}) and
(\ref{3.19-1}), this intersection is not empty if and only~if
$$
\cos\frac{(2\nu+1)\pi}{2j}+\cos\frac{(2l+1)\pi}{2k}=0
$$
for some $\nu\in\{0,\ldots,j-1\}$ and $l\in\{0,\ldots,k-1\},$ which is equivalent to the relation
\begin{equation}\label{4.15}
\frac{2\nu+1}{2j}+(-1)^s\frac{2l+1}{2k}=1+2m
\end{equation}
with $s,\,m\in{\mathbb Z}.$ In its turn, (\ref{4.15}) implies the evenness of $j+k.$ Conversely, let $j+k$ be
even. Then $j$ and $k$ are odd, and (\ref{4.15}) holds for $s=m=0$ with $\nu=(j-1)/2$ and $l=(k-1)/2.$

Finally, let $\alpha=\beta=1.$ Then (\ref{2.8.1}) and (\ref{4.2}) imply that $\det A_{j,k}^{(1,1)}=0$ is
equivalent to ${\cal T}_j\cap\sigma((1/2)A_{1,k}^{(1,1)})\ne\emptyset.$ By virtue of (\ref{3.19}) and
(\ref{3.19-1}), the latter holds if and only if
$$
\cos\frac{(2\nu+1)\pi}{2j}=\cos\frac{l\pi}k
$$
for some $\nu\in\{0,\ldots,j-1\}$ and $l\in\{0,\ldots,k-1\},$ which, in turn, is equivalent to the relation
\begin{equation}\label{4.16}
\frac{2\nu+1}{2j}+(-1)^s\frac{l}k=2m
\end{equation}
for some $s,\,m\in{\mathbb Z}.$ Obviously, (\ref{4.16}) implies the evenness of $k.$ Conversely, let $k$ be
even. Then~$j$ is odd, and (\ref{4.16}) holds for $s=1$ and $m=0$ with $\nu=(j-1)/2$ and $l=k/2.$
$\hfill\Box$
\\

{\large\bf 5. Iso-spectral potentials}
\\

In this section, we return to Inverse Problem~1. The above results give an easy and convenient way for
constructing iso-spectral potentials in the degenerate case.

Let the parameters $\alpha,$ $\beta$ and $a=j/k$ with mutually prime $j$ and $k$ satisfy one of conditions
(I)--(IV) in (\ref{deg}) and, for definiteness, also let $a\in(0,1/2].$ Fix a model complex-valued potential
$q_0(x)\in L(0,1)$ and consider the corresponding eigenvalue problem ${\cal L}(q_0(x),\alpha,\beta,a)$ with
the spectrum $\Lambda:=\{\lambda_n\}_{n\ge1}.$ Consider the set ${\cal M}_\Lambda$ of all corresponding
iso-spectral potentials~$q(x),$ i.e. of such ones for which the spectrum of the problem ${\cal
L}(q(x),\alpha,\beta,a)$ coincides with $\Lambda.$

By virtue of (\ref{2.5}), (\ref{EQ}) and (\ref{2.7}), we have the representation
\begin{equation}\label{5.0}
{\cal M}_\Lambda=\Big\{q_0(x)+g(x): g(x)\in{\cal R}_{j,k}^{(\alpha,\beta)}\Big\},
\end{equation}
where
$$
{\cal R}_{j,k}^{(\alpha,\beta)}=\Big\{R^{-1}F(x): F(x)\in(L(0,b))^k\;\; {\rm and}\;\;
A_{j,k}^{(\alpha,\beta)}F(x)=0 \;\; {\rm a.e.\;\, on}\;\, (0,b)\Big\}, \quad b=\frac1k,
$$
i.e. the supplement $g(x)$ in (\ref{5.0})
is independent of $q_0(x).$

Thus, the question of describing all iso-spectral potentials is reduced to studying the kernel of the matrix
$A_{j,k}^{(\alpha,\beta)}.$ The following lemma answers this question for $j=1.$

\medskip
{\bf Lemma 2. }{\it Each eigenvalue $z_0$ of the matrix $A_{1,k}^{(\alpha,\beta)}$ has the geometric
multiplicity one, while the corresponding eigenvector has the form
\begin{equation}\label{5.1}
X_0=\Big(1,dq_1(z_0),d^2q_2(z_0),\ldots,d^{k-1}q_{k-1}(z_0)\Big)^T.
\end{equation} }

{\it Proof.} According to (\ref{4.6.1}), relation $A_{1,k}^{(\alpha,\beta)}X_0=z_0X_0$ is equivalent to the
system
\begin{equation}\label{5.2}
\left.\begin{array}{l}
[X_0]_1 +d[X_0]_2=z_0[X_0]_1,\\[2mm]
c[X_0]_{m-1} +d[X_0]_{m+1}=z_0[X_0]_m, \quad m=\overline{2,k-1},\\[2mm]
c([X_0]_{k-1} +[X_0]_k)=z_0[X_0]_k.
\end{array}\right\}
\end{equation}
Thus, we have $[X_0]_1\ne0$ as soon as $X_0$ is an eigenvector, otherwise (\ref{5.2}) would imply $X_0=0.$
Without loss of generality, we put $[X_0]_1=1.$ Then the first two lines in (\ref{5.2}) give the relations
\begin{equation}\label{5.3}
[X_0]_2=d(z_0-1), \quad [X_0]_{m+1}=dz_0[X_0]_m-cd[X_0]_{m-1}, \quad m=\overline{2,k-1}.
\end{equation}
Substituting $[X_0]_m=:d^{m-1}Y_m,$ $m=\overline{1,k},$ into (\ref{5.3}), we arrive at the relations
$$
Y_1=1, \quad Y_2=z_0-1, \quad Y_{m+1}=z_0Y_m-cdY_{m-1}, \quad m=\overline{2,k-1}.
$$
Comparing this with (\ref{2.3.3}), we get $Y_m=q_{m-1}(z_0)$ and, hence, $[X_0]_m=d^{m-1}q_{m-1}(z_0)$ for
$m=\overline{1,k},$ which finalizes the proof. $\hfill\Box$

\medskip
It can be easily seen that the last relation in (\ref{5.2}) is fulfilled automatically as soon as $z_0$ is an
eigenvalue of the matrix $A_{1,k}^{(\alpha,\beta)}.$ Indeed, by virtue of (\ref{5.1}), this relation is
equivalent to the relation $cdq_{k-2}(z_0)=(z_0-c)q_{k-1}(z_0),$ which, according to (\ref{2.3.2}), is
equivalent to $p_k(z_0)=0.$

We also note that, according to (\ref{3.6.1}), the algebraic multiplicity of the zero eigenvalue of the
matrix $A_{1,k}^{(0,0)}$ may be equal to $2,$ while, by virtue of Lemma~2, the geometric one cannot.



\medskip
{\bf Lemma 3. }{\it Let the values $\alpha,$ $\beta,$ $j$ and $k$ obey one of conditions~(I)--(IV)
in~(\ref{deg}). Then the kernel of the matrix $A_{j,k}^{(\alpha,\beta)}$ coincides with a linear hull of the
vector $X=(x_1,\ldots, x_k)^T$ determined in the following way:
\begin{equation}\label{5.4}
\;\;\alpha= \beta=0: \quad x_\nu=(-1)^{\nu-1}, \quad \nu=\overline{1,k};
\end{equation}
\begin{equation}\label{5.5}
\alpha=0,\;\; \beta=1: \quad x_\nu=(-1)^{[\frac\nu2]}, \quad \nu=\overline{1,k}; \quad\;
\end{equation}
\begin{equation}\label{5.6}
\alpha=1,\;\; \beta=0: \quad x_\nu=(-1)^{[\frac{\nu-1}2]}, \quad \nu=\overline{1,k}; \;\;\,
\end{equation}
\begin{equation}\label{5.7}
\alpha=
\beta=1: \quad x_\nu=(-1)^{[\frac\nu2]}, \quad \nu=\overline{1,k}.
\end{equation}
}

{\it Proof.} According to Remark~2 in \cite{BK}, in the degenerate case, we have $\rank
A_{j,k}^{(\alpha,\beta)}=k-1,$ i.e. 
$\ker A_{j,k}^{(\alpha,\beta)}$ is always one-dimensional. By virtue of (\ref{2.8.1}) and~(\ref{3.5}), the
matrix $A_{1,k}^{(0,\beta)}$ is degenerate if and only if $\beta=0.$ Hence, relation (\ref{4.1}) along with
Lemma~2 implies that
$\ker A_{j,k}^{(0,0)}$ is a linear hull of the vector $X=X_0$ determined by~(\ref{5.1}) for
$\alpha=\beta=z_0=0.$ Moreover, by virtue of (\ref{4.2}) along with Proposition~1 and Lemma~2, the kernel of
$A_{j,k}^{(1,\beta)}$ for $\beta\in\{0,1\}$ is a linear hull of the vector $X=X_0$ determined by~(\ref{5.1})
for $\alpha=1$ and the corresponding $\beta$ as well as $z_0=0$ since~$j$ is odd in both subcases (III), (IV)
and, hence, $T_j(0)=0.$ Thus, formulae (\ref{5.4}), (\ref{5.6}) and (\ref{5.7}) for components of $X$ can be
easily obtained using (\ref{2.8.1}), (\ref{3.6}) and~(\ref{5.1}).

Now let $(\alpha,\beta)=(0,1).$ Then representation (\ref{4.1}) takes the form
$$
A_{j,k}^{(0,1)}=U_{j-1}\Big(-\frac12A_{1,k}^{(1,0)}\Big)A_{1,k}^{(0,1)}.
$$
According to (II) in (\ref{deg}) as well as (V) in (\ref{non-deg}), we have $\det A_{j,k}^{(0,1)}=0$ if and
only if~$j$ is even. In the degenerate case, since $\det A_{1,k}^{(0,1)}\ne0,$ we have $\det
U_{j-1}((-1/2)A_{1,k}^{(1,0)})=0.$ Moreover, since $j$ and $k$ are mutually prime, the value $k$ is odd.
Thus, by virtue of Proposition~1 along with the relation $U_{j-1}(0)=0,$ a unique up to a multiplicative
constant eigenvector of the matrix~$A_{j,k}^{(0,1)}$ corresponding to the zero eigenvalue satisfies the
linear equation
\begin{equation}\label{5.9}
A_{1,k}^{(0,1)}X=X_0,
\end{equation}
where $X_0$ is an eigenvector of the matrix $A_{1,k}^{(1,0)}$ related to the zero eigenvalue. By virtue of
(\ref{5.6}), we have $[X_0]_\nu=(-1)^{[\frac{\nu-1}2]},$ $\nu=\overline{1,k}.$ Thus, according to
(\ref{4.6.1}), equation (\ref{5.9}) is equivalent to the system of scalar equations
\begin{equation}\label{5.10}
x_1-x_2=1, \quad x_{\nu-1}-x_{\nu+1}=(-1)^{[\frac{\nu-1}2]}, \;\;\nu=\overline{2,k-1}, \quad
x_{k-1}+x_k=(-1)^{[\frac{k-1}2]}.
\end{equation}
Summing up all equations in (\ref{5.10}), we get $2x_1=s_k,$ where
$$
s_n=\sum_{\nu=1}^n(-1)^{[\frac{\nu-1}2]}.
$$
Obviously, $s_{4l+1}=s_{4l+3}=1,$ $s_{4l+2}=2$ and $s_{4l+4}=0$ for all $l\ge0.$ Thus, since $k$ is odd, we
have $s_k=1$ and, hence, $x_1=1/2.$ Rewrite the first $k-1$ equations in (\ref{5.10}) in the following way:
$$
 x_1-x_2=x_1-x_3=1, \quad
 \left.\begin{array}{r} x_{2\nu}-x_{2\nu+2}=(-1)^\nu,\\[2mm] x_{2\nu+1}-x_{2\nu+3}=(-1)^\nu,
 \end{array}\right\}
 \quad\nu=\overline{1,\frac{k-3}2},
$$
whence relations (\ref{5.5}) can be easily established by induction and multiplication with $2.$ $\hfill\Box$

\medskip
{\bf Corollary 5. }{\it Let the values $\alpha,$ $\beta,$ $j$ and $k$ obey one of
conditions~(I)--(IV) in~(\ref{deg}). Then  $F(x)=Xf(x)$ is a general solution of the functional equation
$A_{j,k}^{(\alpha,\beta)}F(x)=0$ in $(L(0,b))^k,$ where components of the vector $X=(x_1,\ldots, x_k)^T$ are
determined by the corresponding formula in (\ref{5.4})--(\ref{5.7}), while the function $f(x)$ ranges over
$L(0,b).$}

\medskip
Thus, we arrive at the following procedure for constructing an iso-spectral potential $q(x)$ that is
different from~$q_0(x).$

\medskip
{\bf Algorithm 1. }{\it Let $q_0(x)\in L(0,1)$ as well as appropriate $\alpha,$ $\beta$ and $j,$ $k$ be
given. Then

(i) Choose a nonzero function $f(x)\in L(0,b);$

(ii) Construct the vector $X=(x_1,\ldots, x_k)^T$ by the corresponding formula in (\ref{5.4})--(\ref{5.7});

(ii) Calculate $q(x)$ by the formulae
\begin{equation}\label{5.11}
q(x)=q_0(x)+R^{-1}F(x), \quad F(x)=Xf(x).
\end{equation}
}

Obviously, the obtained function $q(x)$ ranges over ${\cal M}_\Lambda$ as soon as so does $f(x)$ over
$L(0,b),$ where~$\Lambda$ is the spectrum of the problem ${\cal L}(q_0(x),\alpha,\beta,j/k).$

According to (\ref{2.4.1}) and (\ref{2.3.1}), we have the following formulae for $R^{-1}F(x),$ $x\in(0,1),$
with $F(t)=(f_1(t),\ldots, f_k(t))^T,$ $t\in(0,b):$
$$
R^{-1}F(x)=\left\{\begin{array}{cl} f_\nu(x-(k-\nu)b) \;\; \text{for even} \; j+\nu,\\[3mm]
f_\nu((k-\nu+1)b-x) \;\; \text{for odd} \; j+\nu,
\end{array}\right. \; x\in((k-\nu)b,(k-\nu+1)b), \;\; \nu=\overline{1,k}.
$$
Thus, we arrive at the following representations depending on the parities of~$j$ and~$k:$
\begin{equation}\label{5.12}
R^{-1}F(x)=\left\{\begin{array}{cl} f_k(x), & x\in(0,b),\\
f_{k-1}(2b-x), & x\in(b,2b),\\
f_{k-2}(x-2b), & x\in(2b,3b),\\
f_{k-3}(4b-x), & x\in(3b,4b),\\
f_{k-4}(x-4b), & x\in(4b,5b),\\
\ldots & \\
f_2(1-b-x), & x\in(1-2b,1-b),\\
f_1(x-1+b), & x\in(1-b,1),
\end{array}\right.
\end{equation}
for odd $j$ and odd $k;$
\begin{equation}\label{5.13}
R^{-1}F(x)=\left\{\begin{array}{cl} f_k(b-x), & x\in(0,b),\\
f_{k-1}(x-b), & x\in(b,2b),\\
f_{k-2}(3b-x), & x\in(2b,3b),\\
f_{k-3}(x-3b), & x\in(3b,4b),\\
f_{k-4}(5b-x), & x\in(4b,5b),\\
\ldots & \\
f_2(x-1+2b), & x\in(1-2b,1-b),\\
f_1(1-x), & x\in(1-b,1),
\end{array}\right.
\end{equation}
for even $j$ and odd $k;$
\begin{equation}\label{5.14}
R^{-1}F(x)=\left\{\begin{array}{cl} f_k(b-x), & x\in(0,b),\\
f_{k-1}(x-b), & x\in(b,2b),\\
f_{k-2}(3b-x), & x\in(2b,3b),\\
f_{k-3}(x-3b), & x\in(3b,4b),\\
f_{k-4}(5b-x), & x\in(4b,5b),\\
\ldots & \\
f_3(x-1+3b), & x\in(1-3b,1-2b),\\
f_2(1-b-x), & x\in(1-2b,1-b),\\
f_1(x-1+b), & x\in(1-b,1),
\end{array}\right.
\end{equation}
for odd $j$ and even $k.$
\\

{\large\bf 6. Illustrative examples}
\\

Finally, we give some examples illustrating the term $R^{-1}F(x)$ in (\ref{5.11}) for all degenerate subcases
(I)--(IV) in (\ref{deg}). We also provide the corresponding graphs of $R^{-1}F(x)$ taking the model function
$f(x)$ in (\ref{5.11}) of the following form:
$$
f(x)=\frac{10x}{3b}-\frac{25x^2}{9b^2}, \quad b=\frac1k.
$$

{\bf Example I.} Let $\alpha=\beta=0.$ Then, for $(j,k)=(3,7),$ formulae (\ref{5.4}) and (\ref{5.12}) give
$$
F(x)=\left[\begin{array}{r}1\\-1\\1\\-1\\1\\-1\\1\end{array}\right]f(x),\quad
R^{-1}F(x)=\left\{\begin{array}{rl}f(x), & x\in(0,1/7),\\[2mm]
-f(2/7-x), & x\in(1/7,2/7),\\[2mm]
f(x-2/7), & x\in(2/7,3/7),\\[2mm]
-f(4/7-x), & x\in(3/7,4/7),\\[2mm]
f(x-4/7), & x\in(4/7,5/7),\\[2mm]
-f(6/7-x), & x\in(5/7,6/7),\\[2mm]
f(x-6/7), & x\in(6/7,1),
\end{array}\right.
$$
while, for $(j,k)=(3,8),$ formulae (\ref{5.4}) and (\ref{5.14}) give the representations
$$
F(x)=\left[\begin{array}{r}1\\-1\\1\\-1\\1\\-1\\1\\-1\end{array}\right]f(x),\quad
R^{-1}F(x)=\left\{\begin{array}{rl} -f(1/8-x), & x\in(0,1/8),\\[2mm]
f(x-1/8), & x\in(1/8,1/4),\\[2mm]
-f(3/8-x), & x\in(1/4,3/8),\\[2mm]
f(x-3/8), & x\in(3/8,1/2),\\[2mm]
-f(5/8-x), & x\in(1/2,5/8),\\[2mm]
f(x-5/8), & x\in(5/8,3/4),\\[2mm]
-f(7/8-x), & x\in(3/4,7/8),\\[2mm]
f(x-7/8), & x\in(7/8,1).
\end{array}\right.
$$

\begin{center}
\begin{figure}[!h]
\footnotesize $\qquad$
\begin{minipage}[h]{0.4\linewidth}
\begin{center}
\includegraphics[scale = 0.5]{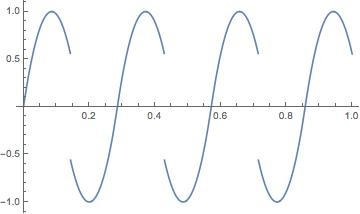} \\[3mm]
Example I: $\alpha=\beta=0,$ $j=3,$ $k=7.$
\end{center}
\end{minipage}
$\qquad\qquad\quad$
\begin{minipage}[h]{0.4\linewidth}
\begin{center}
\includegraphics[scale = 0.5]{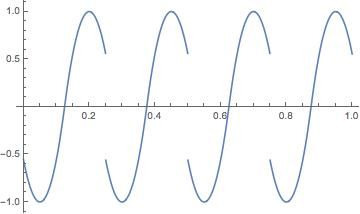} \\[3mm]
Example I: $\alpha=\beta=0,$ $j=3,$ $k=8.$
\end{center}
\end{minipage}
\end{figure}
\end{center}

{\bf Example II.} Let $\alpha=0$ and $\beta=1.$ Then, for $(j,k)=(2,7),$ formulae (\ref{5.5}) and
(\ref{5.13}) give
$$
F(x)=\left[\begin{array}{r}1\\-1\\-1\\1\\1\\-1\\-1\end{array}\right]f(x),\quad
R^{-1}F(x)=\left\{\begin{array}{rl}-f(1/7-x), & x\in(0,1/7),\\[2mm]
-f(x-1/7), & x\in(1/7,2/7),\\[2mm]
f(3/7-x), & x\in(2/7,3/7),\\[2mm]
f(x-3/7), & x\in(3/7,4/7),\\[2mm]
-f(5/7-x), & x\in(4/7,5/7),\\[2mm]
-f(x-5/7), & x\in(5/7,6/7),\\[2mm]
f(1-x), & x\in(6/7,1).
\end{array}\right.
$$

{\bf Example III.} Let $\alpha=1$ and $\beta=0.$ Then, for $(j,k)=(3,7),$ formulae (\ref{5.6}) and
(\ref{5.12}) give
$$
F(x)=\left[\begin{array}{r}1\\1\\-1\\-1\\1\\1\\-1\end{array}\right]f(x),\quad
R^{-1}F(x)=\left\{\begin{array}{rl}-f(x), & x\in(0,1/7),\\[2mm]
f(2/7-x), & x\in(1/7,2/7),\\[2mm]
f(x-2/7), & x\in(2/7,3/7),\\[2mm]
-f(4/7-x), & x\in(3/7,4/7),\\[2mm]
-f(x-4/7), & x\in(4/7,5/7),\\[2mm]
f(6/7-x), & x\in(5/7,6/7),\\[2mm]
f(x-6/7), & x\in(6/7,1).
\end{array}\right.
$$

\begin{center}
\begin{figure}[h]
\footnotesize $\qquad$
\begin{minipage}[h]{0.4\linewidth}
\begin{center}
\includegraphics[scale = 0.5]{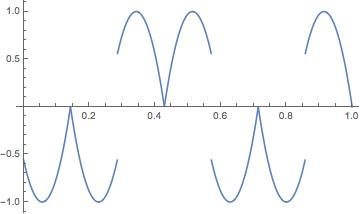} \\[3mm]
Example II: $\alpha=0,$ $\beta=1,$ $j=2,$ $k=7.$
\end{center}
\end{minipage}
$\qquad\qquad\quad$
\begin{minipage}[h]{0.4\linewidth}
\begin{center}
\includegraphics[scale = 0.5]{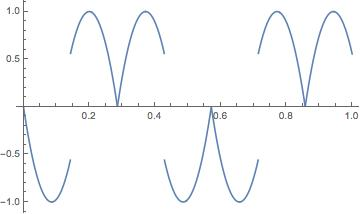} \\[3mm]
Example III: $\alpha=1,$ $\beta=0,$ $j=3,$ $k=7.$
\end{center}
\end{minipage}
\end{figure}
\end{center}

{\bf Example IV.} Let $\alpha=1$ and $\beta=1.$ Then, for $(j,k)=(3,8),$ formulae (\ref{5.7}) and
(\ref{5.14}) give
$$
F(x)=\left[\begin{array}{r}1\\-1\\-1\\1\\1\\-1\\-1\\1\end{array}\right]f(x),\quad
R^{-1}F(x)=\left\{\begin{array}{rl} f(1/8-x), & x\in(0,1/8),\\[2mm]
-f(x-1/8), & x\in(1/8,1/4),\\[2mm]
-f(3/8-x), & x\in(1/4,3/8),\\[2mm]
f(x-3/8), & x\in(3/8,1/2),\\[2mm]
f(5/8-x), & x\in(1/2,5/8),\\[2mm]
-f(x-5/8), & x\in(5/8,3/4),\\[2mm]
-f(7/8-x), & x\in(3/4,7/8),\\[2mm]
f(x-7/8), & x\in(7/8,1).
\end{array}\right.
$$

\begin{center}
\begin{figure}[!h]
\footnotesize $\qquad\qquad\qquad\qquad\qquad\qquad\quad\;\;$
\begin{minipage}[h]{0.4\linewidth}
\begin{center}
\includegraphics[scale = 0.5]{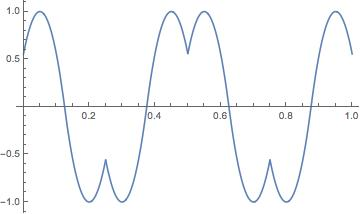} \\[3mm]
$\;$ Example IV: $\alpha=1,$ $\beta=1,$ $j=3,$ $k=8.$
\end{center}
\end{minipage}
\end{figure}
\end{center}

{\bf Acknowledgements.} Sergey Buterin is supported by Grant 20-31-70005 of the Russian Foundation for Basic
Research. Chung-Tsun Shieh is partially supported by the Ministry of Science and Technology, Taiwan under
Grant no. 109-2115-M-032-004-.

\end{document}